\documentclass[11pt,twoside]{article}
\usepackage{amsmath, amssymb, mathrsfs, graphicx, bbm, nicefrac}

\def\az{\alpha}  \def\bz{\beta}

\def\gz{\gamma}  
\def\lz{\lambda}

        \def\uz{\theta}
\def\vz{\varepsilon}

\def\qd{\quad}
\def\qqd{\qquad}

\def\ll{\left}
\def\rr{\right}

\def\sgn{\mathop{\rm sgn}}

\newcommand{\mathsym}[1]{{}}

\def\leq{\leqslant}
\def\geq{\geqslant}

\font\cms=cmss9 scaled \magstep1

\def\nnd{\noindent}

\def\thm{\nnd\bg{thm1}}

\def\lmm{\nnd\bg{lmm1}}
\def\prp{\nnd\bg{prp1}}

\def\xmp{\nnd\bg{xmp1}}

\def\rmk{\nnd\bg{rmk1}}

\def\dethm{\end{thm1}}

\def\delmm{\end{lmm1}}
\def\deprp{\end{prp1}}
\def\dexmp{\end{xmp1}}

\def\dermk{\end{rmk1}}

\def\prf{\medskip \noindent {\bf Proof}. }
\def\deprf{\quad $\square$ \medskip}
\def\bg{\begin}
\def\be{\bg{equation}}
\def\de{\end{equation}}

\def\dear{\end{eqnarray}}
\def\lb{\label}

\newcommand{\rf}[2]{[\ref{#1}; #2]}

\def\den{\end{enumerate}}

\def\d{\text{\rm d}}

\begin{document}

\allowdisplaybreaks[4]
\thispagestyle{empty}
\renewcommand{\thefootnote}{\fnsymbol{footnote}}


\vspace*{.5in}
\begin{center}
{\bf\Large DISCRETE WEIGHTED HARDY INEQUALITIES WITH DIFFERENT KINDS OF BOUNDARY CONDITIONS \footnote{The research is supported by NSFC (Grant No. 11131003) and by the ``985'' project from the Ministry of Education in China.}}
\vskip.15in {Zhong-Wei Liao\footnote{Corresponding author. E-mail: zhwliao@mail.bnu.edu.cn (Tel: +008618810561270)}
}
\end{center}
\begin{center} (School of Mathematical Sciences, Beijing Normal University, Beijing 100875, P. R. China)\\
\vskip.1in
\end{center}


\bigskip

\noindent {\bf Abstract} \qd This paper studies the weighted Hardy inequalities on the discrete intervals with four different kinds of boundary conditions. The main result is the uniform expression of the basic estimate of the optimal constant with the corresponding boundary condition. Firstly, one-side boundary condition is considered, which means that the sequences vanish at the right endpoint (ND-case). Based on the dual method, it can be translated into the case vanishing at left endpoint (DN-case). Secondly, the condition is the case that the sequences vanish at two endpoints (DD-case). The third type of condition is the generality of the mean zero condition (NN-case), which is motivated from probability theory. To deal with the second and the third kinds of inequalities, the splitting technique is presented. Finally, as typical applications, some examples are included. 

\vskip.1in


\vskip.1in


\medskip

\section{Introduction}

In this section, we explain the motivations and main results of this paper. 

Let $[－M, N]$ be a discrete interval $\{-M, -M+ 1, \ldots, N- 1, N\}$ with $M, N\leq \infty$. Here, $[-M, N]$ means $[-M, N)$ if $N= \infty$. In this whole paper, we always assume constants $p, q$ satisfy $1< p\leq q< \infty$ and $\mathbf{u}, \mathbf{v}$ on $[-M, N]$ are two positive sequences. 

We are concerned with the first weighted Hardy inequality
\be\lb{ND Hardy}
\ll[\sum_{n= -M}^N u_n | x_n |^q \rr]^{1/q}\leq A^{ND} \ll[\sum_{n= -M}^N v_n | x_n- x_{n+1} |^p \rr]^{1/p}, \qd x_{N+ 1}=0.
\de
This boundary condition means the sequences vanish at the right endpoint of the interval (i.e. $N+ 1$ is a Dirichlet boundary). When $N= \infty$, $x_{N+ 1}=0$ means that $\lim_{n \rightarrow \infty} x_n =0$, which will not be mentioned again in what follows. The notation $A^{ND}$ is inspired by the relationship between weighted Hardy inequalities and probability theory, cf. \cite{Chen5}. In probabilistic terminology, the corresponding process of this inequality has reflecting (Neumann) boundary at left end point and absorbing (Dirichlet) boundary at right. Naturally, we have the dual version of this case, i.e. the weighted Hardy inequality with Dirichlet boundary at the left endpoint:
\be\lb{DN Hardy}
\ll[\sum_{n= -M}^N u_n | x_n |^q \rr]^{1/q}\leq A^{DN} \ll[\sum_{n= -M}^N v_n | x_n- x_{n-1} |^p \rr]^{1/p}, \qd x_{-M-1}=0.
\de
Besides this one-side condition, we think about the bilateral Dirichlet boundary condition, that is the third weighted Hardy inequality
\be\lb{DD Hardy}
\ll[\sum_{n= -M}^N u_n |x_n|^q \rr]^{1/q}\leq A^{DD} \ll[\sum_{n= -M}^N v_n |x_n- x_{n- 1} |^p \rr]^{1/p}, \qd x_{-M- 1}= x_N= 0.
\de
Finally, we consider the weighted Hardy inequality like
\be\lb{NN Hardy}
\ll[\sum_{n=-M}^N u_n |x_n- m(\mathbf{x})|^q \rr]^{1/q} \leq A^{NN} \ll[\sum_{n=-M}^N v_n |x_n- x_{n- 1}|^p \rr]^{1/p}, 
\de
where $x_{-M- 1}=x_{-M}$ and $m(\mathbf{x})$ is the constant such that 
\be\lb{NN boundary}
\sum_{n= -M}^N u_n \ll| x_n- m(\mathbf{x}) \rr|^{q- 2} \ll(x_n -m(\mathbf{x}) \rr) =0. 
\de
Again, from the probabilistic point of view, this boundary condition corresponds to the ergodic case of processes. 

The main result of this paper is the uniform expression of the basic estimate of the optimal constant $A^{\#}$. For convenience, we give some basic notations. Let $p^*$ be the conjugate number of $p$, i.e. $1/p+ 1/p^* =1$, and $q^*$ is defined similarly. Write $x \wedge y= \min \{x, y\}$ and $x \vee y= \max \{x, y\}$. Define a frequently-used factor 
\be\lb{kqp}
k_{q, p}= \ll[\frac{q- p}{p B\ll(\frac{p}{q-p}, \frac{p(q-1)}{q-p} \rr)}\rr]^{1/p- 1/q}, q>p>1; \qd k_{p, p}= p^{1/p} (p^*)^{1/p^*}, p>1, 
\de
where $B(a, b)= \int_0^1 x^{a- 1} (1- x)^{b- 1} \d x$ is the Beta function. In this whole paper, we adopt some usual conventions. Firstly, the absolute value function $f(x)= |x|$ is  differentiable everywhere except for $x=0$, and then we define the derivative formally 
$$
\ll(|x| \rr)' = \sgn (x). 
$$
Secondly, the definition of $x_{\infty}$ is
$$
x_{\infty}:= \lim_{n \rightarrow \infty} x_n.
$$
Finally, we define the summation on a empty set is zero, i.e. 
$$
\sum_{i\in \emptyset} x_i =0.  
$$
Here is the main theorem. 

\thm\lb{main thm}
{\cms
Write $\hat{v}_i= v_i^{1-p^*}$ ($-M \leq i\leq N$). For the optimal constant $A^{\#}$ in the weighted Hardy inequalities with four different boundary conditions, we have the following uniform expression of the basic estimate 
$$
B_*^{\#}\leq A^{\#} \leq k_{q, p} B^{* \#}, 
$$
where $k_{q, p}$ is defined as (\ref{kqp}) and 
$$
\aligned
B_*^{ND} & = B^{* ND} = \sup_{n\in [-M, N]} \ll(\sum_{i= -M}^n u_i \rr)^{1/q} \ll(\sum_{j= n}^N \hat{v}_j \rr)^{1/p^*}, \\
B_*^{DN} & = B^{* DN} = \sup_{n\in [-M, N]} \ll(\sum_{i= -M}^n \hat{v}_i \rr)^{1/p^*} \ll(\sum_{j=n}^N u_j \rr)^{1/q}, \\
B_*^{DD} & = \sup_{-M\leq x< y\leq N} \ll(\sum_{n=x}^{y- 1} u_n \rr)^{1/q} \Bigg[\Bigg(\sum_{i= -M}^x \hat{v}_i \Bigg)^{1- p} + \Bigg(\sum_{j= y}^N \hat{v}_j \Bigg)^{1- p} \Bigg]^{-1/p}, \\
B^{* DD} & = \sup_{-M\leq x< y\leq N} \ll(\sum_{n=x}^{y- 1} u_n \rr)^{1/q} \Bigg[\Bigg(\sum_{i= -M}^x \hat{v}_i \Bigg)^{-q/p^*} + \Bigg(\sum_{j= y}^N \hat{v}_j \Bigg)^{-q/p^*} \Bigg]^{-1/q}, \\
B_*^{NN} & = \sup_{-M\leq x< y\leq N} \ll(\sum_{n=x+1}^{y} \hat{v}_n \rr)^{1/p^*} \Bigg[\Bigg(\sum_{i= -M}^x u_i \Bigg)^{1- q^*} + \Bigg(\sum_{j= y}^N u_j \Bigg)^{1- q^*} \Bigg]^{-1/q^*}, \\
B^{* NN} &= \sup_{-M\leq x< y\leq N} \ll(\sum_{n=x+ 1}^{y} \hat{v}_n \rr)^{1/p^*} \Bigg[\Bigg(\sum_{i= -M}^x u_i \Bigg)^{-p^*/q} + \Bigg(\sum_{j= y}^N u_j \Bigg)^{-p^*/q} \Bigg]^{-1/p^*}.  
\endaligned
$$
}
\dethm

\rmk\lb{main rmk}
{\cms

(1) The result of DN-case can refer to \rf{Liao}{Corollary 1.2 and Theorem 1.4}, and the factor $k_{q, p}$ is sharp when $N= \infty$ and $\sum_{i= -M}^\infty u_i = \infty$. 

(2) Correspondingly, in the ND-case, the factor $k_{q, p}$ is sharp when $N= \infty$ and $\sum_{i= -M}^\infty \hat{v}_i = \infty$. 

(3) In the DD-case and the NN-case, the upper estimates effective when $1< p\leq q< \infty$. However, we need $1< p, q< \infty$ only for the lower estimates. 

(4) For DD-case and NN-case, when $p\leq q$, we have 
$$
B_*^{\#} \leq B^{* \#} \leq 2^{1/p- 1/q} B_*^{\#}. \\
$$
In particular, when $p= q$, we have $B_*^{DD} =B^{* DD}$ and $B_*^{NN} = B^{* NN}$. 
}
\dermk

In the following part, we will give the motivations and the advance of the related studies.

The first thing to notice here is the the duality between the ND-case and the DN-case. About this pair of problems, there are lots relevant results in recent research, cf. \cite{Chen5, Chen3, Chen2, Wang, Maz'ya, Opic}. For the continuous version of the DN-case (\ref{DN Hardy}), B. Opic and A. Kunfner \rf{Opic}{Theorem 1.14} and V. G. Maz'ya \rf{Maz'ya}{Theorem 1, pp. 42-43} gave the well know basic estimate
\be\lb{basic estimate}
B\leq A \leq \tilde{k}_{q, p} B,
\de
where the factor $\tilde{k}_{q, p}$ is constant
$$
\tilde{k}_{q, p}= \ll(1+ \frac{q}{p^*} \rr)^{1/q} \ll(1+ \frac{p^*}{q} \rr)^{1/p^*},
$$
which is bigger than $k_{q, p}$. Afterwards, the variational formula of the optimal constants was discovered by Chen \rf{Chen3}{Theorem 2.1}, which is an entirely different method. As a direct application of this method, Chen gave the basic estimate \rf{Chen3}{Corollary 2.3}, which is consistent with (\ref{basic estimate}), and the approximating procedure\rf{Chen3}{Theorem 2.2}. It's worth mentioning that, with the result of Bliss \cite{Bliss}, the factor $\tilde{k}_{q,p}$ in (\ref{basic estimate}) can be improved to $k_{q, p}$, cf. \rf{Bennett3}{Theorem 8}, \rf{Manakov}{Theorem 2} and \rf{Kufner}{pp. 45-47}. The above mentioned results are also true for discrete case, which are main work of \cite{Liao}. For the continuous version of the ND-case (\ref{ND Hardy}), B. Opic and A. Kunfner \rf{Opic}{Remark 1.8} gave a straightforward transformation so that we can reduce the ND-case to the DN-case, that means there are corresponding conclusions of inequality (\ref{ND Hardy}), cf. \rf{Opic}{Theorem 6.2}. Based on the dual relationship of these two inequalities, B. Opic and A. Kunfner \rf{Opic}{Chapter 1, Section 7} gave a series of results of this pair of inequalities by dual method. As an application of this method in discrete case, we will give the the proof of the ND-case in section 2. 

Attention will turn next to the bilateral Dirichlet boundary weighted Hardy inequalities. In the special case that $p= q= 2$, there are lots of results are proved by Chen, cf. \cite{Chen5} (for discrete case) and \cite{Chen9} (for continuous case). For the general case $1< p\leq q< \infty$, the basic estimate of the optimal constant in the continuous case is given in \cite{Chen4}. In addition, P. Gurka \cite{Gurka} (also see \rf{Opic}{Theorem 8.2}) shown that this bilateral Dirichlet boundary weighted Hardy inequalities holds if and only if 
\be\lb{B_Opic}
B:= \sup_{-M \leq x< y\leq N} \ll(\sum_{n=x}^{y- 1} u_n \rr)^{1/q} \Bigg[\Bigg(\sum_{i= -M}^x \hat{v}_i \Bigg)^{1/p^*} \wedge \Bigg(\sum_{j= y}^N \hat{v}_j \Bigg)^{1/p^*} \Bigg] < \infty.
\de
Based on the result, B. Opic and A. Kunfner \rf{Opic}{Remark 8.3} gave the following estimate
$$
2^{-1/p} B \leq A^{DD} \leq 2 (2^q - 1)^{1/q} B. 
$$
We can compare this result with Theorem \ref{main thm} by direct calculation. Using the following inequalities 
$$
\az \vee \bz \leq \az + \bz \leq 2 (\az \vee \bz), \qd \text{for all $\az > 0$, $\bz> 0$}, 
$$
$$
k_{q, p} \leq \tilde{k}_{q, p} \leq k_{q, q} \leq 2 (2^q - 1)^{1/q}, \qd \text{$1< p\leq q< \infty$}, 
$$
we obtain 
$$
2^{-1/p} B \leq B_*^{DD} \leq A^{DD} \leq k_{q, p} B^{* DD} \leq 2 (2^q - 1)^{1/q} B.
$$
The main method of the DD-case of Theorem \ref{main thm} is the splitting technique. The idea is that the space can be split into two parts so that the problem with bilateral vanishing boundaries degenerate into one-side boundary, cf. \cite{Chen6, Chen4, Chen7, Chen8, Mao2}. The construction of this part is inspired by the results proved by Chen, Zhang and Zhao \cite{Chen8}. The details will be provided in the Section 3. 

Finally, the main concern is the NN-case. The additional condition (\ref{NN boundary}), which is the generality of the mean zero condition. For detail, we assume that $q=2$ and $\sum_{i= -M}^N u_i < \infty$, then we can define a discrete probability measure $\mathbf{\pi}$ as
$$
\pi_n= \frac{u_n}{\sum_{i= -M}^N u_i}. 
$$
It's easy to see that $m(\mathbf{x})$ is the mean value of $\mathbf{x}$ with respect to $\mathbf{\pi}$, and sequence $\mathbf{x}- m(\mathbf{x})$ satisfies the mean zero condition with respect to $\mathbf{u}$. In the special case $p= q$, L. D. Wang had worked some related results, which had not been officially published but gave the author much guidance in the research of $p \leq q$ case. From another point of view, the DD-case (\ref{DD Hardy}) and the NN-case (\ref{NN Hardy}) inequalities have its own spectral meaning. Write $x'_n = x_{n}- x_{n-1}$. When $p= q$, $A^{-p}$, which is the $-p$ power of the optimal constant, in (\ref{DD Hardy}) (resp., (\ref{NN Hardy})) corresponds to the first nontrivial eigenvalue $\lz$ of 
\be\lb{characteristic equation}
(-1) \lz u_n |x_n|^{q- 1} \sgn(x_n) = v_{n+ 1} |x'_{n+ 1}|^{p- 1} \sgn(x'_{n+ 1}) - v_{n} |x'_n|^{p- 1} \sgn(x'_n)
\de
with boundary condition $x_{-M- 1}= x_N= 0$ (resp., $x'_{-M}= x'_{N+ 1}= 0$). For more details about this view, readers can refer to the \cite{Buttazzo}, \cite{Gurka} and the series researches of Chen, cf. \cite{Chen1, Chen6, Chen5, Chen9, Chen7, Chen8}. In a similar way, we use the splitting technique to prove the basic estimate of this optimal constant. 

\section{Proof of the ND-case}

As previously mentioned, we use the dual method to prove the first part of Theorem \ref{main thm}. For any non-negative sequence $\mathbf{x}$ on $[-M, N]$, define
\be\lb{H}
H \mathbf{x} (n) = \sum_{i=-M}^n x_i, \ \ n\geq -M; \qd H^* \mathbf{x} (n) = \sum_{i=n}^N x_i, \ \ n\leq N.
\de
Obviously, $H^*$ is the adjoint operator of $H$ with respect to the bilinear function $\langle \cdot, \cdot \rangle$, which is defined as  
\be\lb{<,>}
\langle \mathbf{x}, \mathbf{y} \rangle= \sum_{n=-M}^N x_n y_n, 
\de
where $\mathbf{x}$ and $\mathbf{y}$ are arbitrary non-negative sequences on $[-M, N]$.

\prp\lb{Decreasing}
{\cms
Given a sequence $\mathbf{x}$ on $[-M, N]$, define a decreasing sequence $y_n = \sup_{k\geq n} \ll\{|x_k| \rr\}$, then we have 
$$
\frac{\ll[\sum_{n=-M}^N u_n |x_n|^q \rr]^{1/q}}{\ll[\sum_{n= -M}^N v_n |x_n- x_{n+1}|^p\rr]^{1/p}} \leq \frac{\ll[\sum_{n=-M}^N u_n |y_n|^q \rr]^{1/q}}{\ll[\sum_{n=-M}^N v_n |y_n - y_{n+1}|^p \rr]^{1/p}}. 
$$
}
\deprp

\prf
The proof is easy by direct calculation. 
\deprf

Based on this proposition, it would be sufficient to think about the decreasing sequences when we consider the basic estimate of the optimal constant of (\ref{ND Hardy}). In the same way, we only need to think about the increasing sequences of (\ref{DN Hardy}). Through this result and the aforementioned symbols, we can rewrite (\ref{ND Hardy}) and (\ref{DN Hardy}) as 
\be\lb{ND Hardy2}
\ll[\sum_{n= -M}^N u_n \ll| H^* \mathbf{x} (n) \rr|^q \rr]^{1/q} \leq A^{ND} \ll[\sum_{n= -M}^N v_n \ll|x_n \rr|^p \rr]^{1/p}
\de
and 
\be\lb{DN Hardy2}
\ll[\sum_{n= -M}^N u_n \ll| H \mathbf{x} (n) \rr|^q \rr]^{1/q} \leq A^{DN} \ll[\sum_{n= -M}^N v_n \ll|x_n \rr|^p \rr]^{1/p}. 
\de

As mentioned in Remark \ref{main rmk}, about the DN-case, we already have the variational formulas of the optimal constants, the approximating procedure and the basic estimates with the improved factor, cf. \cite{Liao}. For convenience, we reiterate the DN-case of Theorem \ref{main thm} as the following lemma. For details, see \rf{Liao}{Corollary 1.2 and Theorem 1.4}. 

\lmm\lb{DN lmm}
{\cms
The inequality (\ref{DN Hardy2}) holds for every sequence $\mathbf{x}$ if and only if $B^{DN}< \infty$, where 
\be\lb{DN B}
B^{DN}= \sup_{n\in [-M, N]} \ll(\sum_{i= -M}^n \hat{v}_i \rr)^{1/p^*} \ll(\sum_{j=n}^N u_j \rr)^{1/q}.
\de
Moreover, we have
\be\lb{DN estimate}
B^{DN} \leq A^{DN} \leq k_{q, p} B^{DN}, 
\de
where $k_{q, p}$ is defined as (\ref{kqp}). In particular, the factor $k_{q, p}$ is sharp when $N= \infty$ and $\sum_{i= -M}^{\infty} \hat{v}_i = \infty$.
}
\delmm

Define 
\be\lb{l^q space}
l^q (\mathbf{u})= \ll\{\text{ $\mathbf{x}$ on $[-M, N]$}: \text{$x_n\geq 0$ and } \|\mathbf{x} \|_{l^q(\mathbf{u})}^q := \sum_{n= -M}^N u_n x_n^q < \infty \rr\}, 
\de
Obviously, with norm $\| \cdot \|_{l^q(\mathbf{u})}$, $l^q (\mathbf{u})$ becomes a Banach space. Similarly, we can define a Banach space $l^p(\mathbf{v})$ with norm $\| \cdot \|_{l^p(\mathbf{v})}$. Write $\mathbf{u}^{1- q^*}= (u_{-M}^{1- q^*}, \dots, u_{N}^{1- q^*})$. By direct calculation, it is easy to show that $l^{q^*}(\mathbf{u}^{1- q^*})$ (resp., $l^{p^*}(\mathbf{v}^{1- p^*})$) is the conjugate space of $l^{q}(\mathbf{u})$ with respect to the inner product $\langle \cdot, \cdot \rangle$.

\medskip
\noindent {\bf Proof of the ND-case.} Since $H^*$ is the adjoint operator of $H$, we have
$$
A^{ND} = \|H^* \|_{l^p(\mathbf{v}) \rightarrow l^q(\mathbf{u})} = \|H \|_{l^{q^*}(\mathbf{u}^{1- q^*}) \rightarrow l^{p^*}(\mathbf{v}^{1- p^*})}.
$$
Hence, inequality (\ref{ND Hardy2}) holds if and only if the following inequality holds: 
\be\lb{= 2_3}
\ll[\sum_{n=-M}^N v_n^{1- p^*} \ll(H \mathbf{x} (n) \rr)^{p^*} \rr]^{1/p^*} \leq A^{ND} \ll[\sum_{n=-M}^N u_n^{1- q^*} x_n^{q^*} \rr]^{1/q^*}. 
\de
Since $1< p\leq q< \infty$, we have $1< q^*\leq p^*< \infty$. As an application of Lemma \ref{DN lmm}, the inequality (\ref{= 2_3}) holds if and only if 
\be\lb{= 2_4}
\sup_{n\in [-M, N]} \ll(\sum_{i=-M}^n  (u_i^{1- q^*})^{1- q} \rr)^{1/q} \ll(\sum_{j= n}^N v_j^{1-p^*} \rr)^{1/{p^*}} < \infty.
\de
Hence, define $B^{ND}$ by (\ref{= 2_4}) and we get the conclusion immediately.
\deprf

\section{The Proof of the DD-case}

It all started with splitting technique. The main idea of the following construction is inspired by the result of \cite{Chen8}.

Given any $\uz \in (-M, N)$, we can construct two inequalities on the left- and the right-hand sides of $\uz$. Fix a constant $0\leq \gz \leq 1$, define $\mathbf{u}^-$ on $[-M, \uz]$ and $\mathbf{u}^+$ on $[\uz+1, N+ 1]$ as
$$
{u}^-_n =
\begin{cases}
u_n, & n\in [-M, \uz- 1], \\
(1- \gz) u_{\uz}, & n= \uz,
\end{cases}
{u}^+_n =
\begin{cases} 
\gz u_{\uz}, \qqd & n= \uz+ 1, \\
u_{n- 1}, \qqd & n\in [ \uz+2,  N+1].
\end{cases}
$$
The following lemma is the key in splitting technique. 

\lmm\lb{calculation}
{\cms
Given a non-negative sequence $\mathbf{x}$ on $[-M, N]$, satisfying $x_{-M-1}= x_N= 0$. Define $\mathbf{x}^-$ on $[-M, \uz]$ and $\mathbf{x}^+$ on $[\uz+ 1, N+ 1]$ as follows
\be\lb{x}
x^-_n = x_n, -M\leq n\leq \uz, \qqd x^+_n = x_{n- 1}, \uz+ 1\leq n\leq N+1.
\de
Then we have 
$$
\aligned
\sum_{n=-M}^N &u_n |x_n|^q = \sum_{n= -M}^{\uz} {u}^-_n |x^-_n|^q+ \sum_{n= \uz+ 1}^{N+ 1} u^+_n |x^+_n|^q; \\
\sum_{n=-M}^N &v_n |x_n- x_{n- 1}|^p = \sum_{n=-M}^{\uz} v_n |x^-_n- x^-_{n- 1}|^p + \sum_{n= \uz+ 1}^{N} v_n |x^+_n- x^+_{n+ 1}|^p.
\endaligned
$$
}
\delmm

Based on this result, we can construct two inequalities
\begin{align}
\ll[\sum_{n= -M}^{\uz} u^-_n |x_n|^q \rr]^{1/q}\leq &A^-(\uz, \gz) \ll[\sum_{n= -M}^{\uz} v_n |x_n- x_{n- 1} |^p \rr]^{1/p}, \  x_{-M- 1}= 0, \lb{DD A^-}\\
\ll[\sum_{n= \uz+ 1}^{N} u^+_n |x_n|^q \rr]^{1/q}\leq &A^+(\uz, \gz) \ll[\sum_{n= \uz+ 1}^N v_n |x_n- x_{n+ 1} |^p \rr]^{1/p}, \  x_{N+ 1}= 0, \lb{DD A^+}
\end{align}
where $A^-(\uz, \gz)$ and $A^+(\uz, \gz)$ are recorded as the optimal constants of the corresponding inequalities. Clearly, these two inequalities only have one-side boundary conditions. The first result is about the relationship between the optimal constant $A$ in (\ref{DD Hardy}) and the optimal constants $A^\pm(\uz, \gz)$ of the inequalities (\ref{DD A^-}) and (\ref{DD A^+}). 

\prp\lb{relationship}
{\cms
Let $1< p\leq q< \infty$, we have
$$
\aligned
2^{1/q- 1/p} &\sup_{\uz \in [-M, N]} \sup_{0\leq \gz \leq 1} \ll(A^-(\uz, \gz) \wedge A^+(\uz, \gz) \rr) \leq A \\
& \leq \inf_{\uz \in [-M, N]} \inf_{0\leq \gz \leq 1} \ll(A^-(\uz, \gz) \vee A^+(\uz, \gz) \rr).
\endaligned
$$
}
\deprp

\prf
(a) For any $\uz\in [-M, N]$, $\gz\in [0, 1]$ and any non-negative sequence $\mathbf{x}$ on $[-M, N]$ satisfying $x_{-M- 1}= x_N= 0$, by Lemma \ref{calculation}, (\ref{DD A^-}) and (\ref{DD A^+}), we have 
$$
\aligned
\sum_{n=-M}^N &v_n |x_n- x_{n- 1}|^p = \sum_{n=-M}^{\uz} v_n |x^-_n- x^-_{n- 1}|^p + \sum_{n= \uz+ 1}^{N} v_n |x^+_n- x^+_{n+ 1}|^p \\
&\geq \ll( A^-(\uz, \gz) \rr)^{-p} \ll(\sum_{n= -M}^\uz u^-_n |x^-_n|^q \rr)^{p/q} + \ll( A^+(\uz, \gz) \rr)^{-p} \ll(\sum_{n= \uz+ 1}^N u^+_n |x^+_n|^q \rr)^{p/q} \\
&\geq \ll[\ll( A^-(\uz, \gz) \rr)^{-p} \wedge \ll( A^+(\uz, \gz) \rr)^{-p} \rr] \ll(2^{(p/q -1) \vee 0} \rr)^{-1} \ll[\sum_{n= -M}^{N} u_n |x_n|^q \rr]^{p/q}.
\endaligned
$$
The last step is $c_r$-inequality. Since $\mathbf{x}$, $\uz$ and $\gz$ are arbitrary, we obtain
$$
A \leq 2^{(1/q- 1/p) \vee 0} \inf_{\uz\in [-M, N]} \inf_{0\leq \gz \leq 1} \ll( A^-(\uz, \gz) \vee A^+(\uz, \gz) \rr). 
$$
Clearly, $2^{(1/q- 1/p) \vee 0}= 1$ when $1< p\leq q< \infty$, which competes the upper bounds of $A$. 

(b) For any $\vz> 0$, since $A^-(\uz, \gz)$ is optimal constant, we can construct a sequence $\mathbf{x^-}$ satisfies: 
$$
\aligned
& \sum_{n=-M}^{\uz} u^-_n |x^-_n|^q =1, \qd x^-_{-M- 1}= 0\\
&\ll[\sum_{n= -M}^\uz v_n |x^-_n- x^-_{n- 1}|^p \rr]^{1/p} < \ll( A^-(\uz, \gz) \rr)^{-1} + \vz. 
\endaligned
$$
Moreover, we assert $x^-_\uz > 0$. If not, i.e. $x^-_\uz = 0$, we can construct another sequence $\mathbf{\bar{x}}$ as
$$
\bar{x}_n= 
\begin{cases}
x^-_n /m, & -M-1 \leq n \leq a, \\
x^-_a /m, & a< n \leq \uz, 
\end{cases}
$$
where $a\in (-M-1, \uz)$ satisfies $x^-_a= \displaystyle \max_{-M-1< k< \uz} \{x^-_k\} >0$ and $m$ is a constant defined as
$$
m := \ll[ \sum_{n= -M-1}^{a} u_n^- |x_n^-|^q + \sum_{n= a+1}^{\uz} u_n^- |x_a^-|^q \rr] ^{1/q} \geq 1. 
$$
Then we have
$$
\sum_{n=-M}^{\uz} u^-_n |\bar{x}_n|^q =1,
$$
$$
\ll[\sum_{n= -M}^\uz v_n |\bar{x}_n- \bar{x}_{n- 1}|^p \rr]^{1/p} = \frac{1}{m} \ll[\sum_{n= -M}^a v_n |x_n^- - x_{n- 1}^- |^p \rr]^{1/p} < \ll( A^-(\uz, \gz) \rr)^{-1} + \vz. 
$$
Therefore, we can use $\mathbf{\bar{x}}$ instead of $\mathbf{x^-}$ when $x^-_\uz= 0$. 

Use the same method, we can construct a sequence $\mathbf{x}^+$ satisfies
$$
\aligned
&\sum_{n=\uz+ 1}^{N} u^+_n |x^+_n|^q =1, \qd x^+_{N+ 1}= 0, \qd x^+_{\uz+ 1}> 0, \\
&\ll[\sum_{n=\uz+ 1}^{N} v_n |x^+_n- x^+_{n+ 1}|^p \rr]^{1/p} < \ll( A^+(\uz, \gz) \rr)^{-1} + \vz. 
\endaligned
$$
From sequences $\mathbf{x^+}$ and $\mathbf{x^-}$, we can obtain the needed conclusion. In detail, set $c= x^+_{\uz+ 1}/ x^-_{\uz}$, define $\mathbf{x}$ on $[-M, N]$ as
$$
x_n= 
\begin{cases}
c x^-_n, & -M\leq n\leq \uz, \\
x^+_{n+1}, & \uz+ 1\leq n\leq N.
\end{cases}
$$
By Lemma \ref{calculation} we obtain
\be\lb{**}
\sum_{n=-M}^N u_n |x_n|^q = c^q \sum_{n=-M}^\uz u^-_n |x^-_n|^q + \sum_{n= \uz+1}^{N} u^+_n |x^+_n|^q = c^q+ 1, 
\de
and 
$$
\aligned
\Big( \sum_{n=-M}^N v_n |x_n & - x_{n-1}|^p \Big)^{\frac{1}{p}} = \Big( c^p \sum_{n=-M}^\uz v_n |x^-_n- x^-_{n-1}|^p + \sum_{n=\uz+ 1}^N v_n |x^+_n- x^+_{n+1}|^p \Big)^{\frac{1}{p}} \\
&\leq \ll(c^p+ 1 \rr)^{1/p} \ll[\ll( A^-(\uz, \gz) \rr)^{-1} \vee \ll( A^+(\uz, \gz) \rr)^{-1} + \vz \rr] \\
&\leq 2^{1/p- 1/q} \ll(c^q+ 1 \rr)^{1/q} \ll[\ll( A^-(\uz, \gz) \rr)^{-1} \vee \ll( A^+(\uz, \gz) \rr)^{-1} + \vz \rr]. 
\endaligned
$$
Combine (\ref{**}), since $\uz$ and $\gz$ are arbitrary, let $\vz \rightarrow 0$ we have
$$
A \geq 2^{1/q- 1/p} \sup_{\uz \in [-M, N]} \sup_{0\leq \gz \leq 1} \ll[A^-(\uz, \gz) \wedge A^+(\uz, \gz) \rr]. 
$$
This completes the proof of Proposition \ref{relationship}.
\deprf

As applications of the basic estimates of the optimal constant in the ND-case and the DN-case of Theorem \ref{main thm}, we have
\be\lb{DD estimate}
B^{\pm}(\uz, \gz) \leq A^{\pm}(\uz, \gz) \leq k_{q, p} B^{\pm}(\uz, \gz), 
\de
where $k_{q, p}$ is defined as (\ref{kqp}) and $B^{\pm}(\uz, \gz)$ are given below: 
\begin{align}\lb{B^+-}
B^-(\uz, \gz) &= \sup_{n\in [-M, \uz]} \ll(\sum_{i=-M}^n \hat{v}_i \rr)^{1/p^*} \ll(\sum_{j=n}^{\uz-1} u_j + (1- \gz) u_{\uz}\rr)^{1/q}, \\
B^+(\uz, \gz) &= \sup_{n\in [\uz+ 1, N]} \ll(\sum_{i= \uz+ 1}^{n- 1} u_i + \gz u_{\uz}\rr)^{1/q} \ll(\sum_{j=n}^N \hat{v}_j \rr)^{1/p^*}.
\end{align}

\medskip
\noindent {\bf Proof of the DD-case.} For brevity, we use $B^*$ (resp., $B_*$) instead of $B^{* DD}$ (resp., $B_*^{DD}$) in this proof.  

(a) By proportional property and direct calculation, we give the relationship between $B^*$ and $B^\pm (\uz, \gz)$
\be\lb{DD =1}
\sup_{\uz \in[-M, N]} \sup_{0\leq \gz \leq 1} \ll[B^-(\uz, \gz) \wedge B^+(\uz, \gz) \rr] \leq B^*. 
\de

(b) Next, we assert that exist $\bar{\uz}$ and $\bar{\gz}$ such that 
\be\lb{DD =2}
B^-(\bar{\uz}, \bar{\gz})= B^+(\bar{\uz}, \bar{\gz}). 
\de
In fact, fix $0\leq \gz \leq 1$, it is easy to prove $B^-(\uz, \gz)$ is increasing with respect to $\uz$, and $B^+(\uz, \gz)$ is decreasing. Besides, with direct calculation, we have
\be\lb{DD =3}
B^-(\uz, 0)= B^-(\uz+1, 1), \qqd B^+(\uz, 0)= B^+(\uz+1, 1). 
\de
Fix $\gz=0$, since the monotonicity of $B^-(\uz, \gz)$ and $B^+(\uz, \gz)$, there is a point satisfies  
$$
\bar{\uz}= \max \{-M- 1\leq n\leq N- 1: B^+(n, 0)\geq B^-(n, 0)\}. 
$$ 
Then, by (\ref{DD =3}), we have
\be\lb{DD =4}
\begin{cases}
B^+(\bar{\uz}, 0)\leq B^-(\bar{\uz}, 0), \\
B^+(\bar{\uz}, 1)\geq B^-(\bar{\uz}, 1).
\end{cases}
\de
It's important to note that, for any given $\bar{\uz}$, $B^-(\bar{\uz}, \gz)$ continuously decrease respect to $\gz$, and $B^+(\bar{\uz}, \gz)$ continuously increase. Hence, by (\ref{DD =4}), the existence of $\bar{\gz}$ is obvious.

Combining (\ref{DD =1}) and (\ref{DD =2}), we have 
\be\lb{a+b}
B^-(\bar{\uz}, \bar{\gz}) \leq \sup_{\uz \in[-M, N]} \sup_{0\leq \gz \leq 1} \ll[B^-(\uz, \gz) \wedge B^+(\uz, \gz) \rr] \leq B^*.
\de

(c) Using Proposition \ref{relationship}, the ND-case and DN-case of Theorem \ref{main thm}, Lemma \ref{DN lmm} and (\ref{a+b}), we obtain the upper bound of $A$ immediately
$$
\aligned
A &\leq \inf_{\uz\in [-M, N]} \inf_{0\leq \gz\leq 1} \ll(A^-(\uz, \gz) \vee A^+(\uz,\gz) \rr) \\
&\leq k_{q,p} \inf_{\uz\in [-M, N]} \inf_{0\leq \gz\leq 1} \ll(B^-(\uz, \gz) \vee B^+(\uz,\gz) \rr) \\
&\leq k_{q,p} B^-(\bar{\uz}, \bar{\gz}) \leq k_{q,p} B^*.
\endaligned
$$

(d) For the lower bound of $A$, we use a straightforward method to get the conclusion. Given $x$, $y$ with $-M\leq x< y\leq N$, define $\mathbf{\tilde{x}}$ on $[-M, N+1]$ as 
\be\lb{DD =5}
\tilde{x}_n= 
\begin{cases}
c \sum_{i= -M}^n \hat{v}_i, & -M\leq n\leq x, \\
\sum_{i=y}^N \hat{v}_i, & x< n< y+ 1, \\
\sum_{i=n}^N \hat{v}_i, & y+ 1\leq n\leq N+1, 
\end{cases}
\de
where $c= \ll(\sum_{i=y}^N \hat{v}_i \rr)/ \ll(\sum_{i=-M}^x \hat{v}_i \rr)$. Clearly, we have the boundary conditions $\tilde{x}_{-M-1}= 0$ and $\tilde{x}_{N+ 1}= 0$. For any $\uz \in [x, y]$, define $\mathbf{x}$ on $[-M, N]$ as
$$
x_n= 
\begin{cases}
\tilde{x}_n,  & -M\leq n\leq \uz, \\
\tilde{x}_{n+ 1}, & \uz+ 1 \leq n \leq N, 
\end{cases}
$$
then applying Lemma \ref{calculation}, we have
$$
\aligned
\sum_{k= -M}^N u_k |x_k|^q &\geq \sum_{k= x}^y u_k |x_k|^q = \ll(\sum_{k=x}^y u_k \rr) \ll(\sum_{i= y}^N \hat{v}_k\rr)^q, \\
\sum_{k=-M}^N v_k |x_k- x_{k-1}|^p &= c^p \ll(\sum_{k=-M}^{x} \hat{v}_k\rr) + \ll(\sum_{k=y}^{N} \hat{v}_k\rr). 
\endaligned
$$
Hence, 
$$
A\geq \ll(\sum_{k=x}^y u_k \rr)^{1/q} \ll[\ll(\sum_{k=-M}^{x} \hat{v}_k\rr)^{1-p} + \ll(\sum_{k=y}^{N} \hat{v}_k\rr)^{1-p}\rr]^{-1/p}. 
$$
Since $x$ and $y$ are arbitrary, we obtain $A \geq B_*$.

(e) In the remaining part of this proof, we consider the relationship of $B^*$ and $B_*$. Obviously, we have $B_*=B^*$ when $p= q$. In the case of $1< p\leq q< \infty$, we use $c_r$ inequality: 
\be\lb{DD =6}
\ll(\az+ \bz\rr)^r \leq 2^{(r- 1)\vee 0} \ll(\az^r+ \bz^r\rr), \qqd \az, \bz\in \mathbb{R}^+.
\de
Set 
$$
\az= \ll(\sum_{i=-M}^x \hat{v}_i \rr)^{-q/p*}, \qd \bz= \ll(\sum_{i=y}^N \hat{v}_i \rr)^{-q/p*}, \qd r=\frac{p}{q} \leq 1, 
$$
then we get $B_*\leq B^*$.
On the other hand, set 
$$
\az= \ll(\sum_{i=-M}^x \hat{v}_i \rr)^{1- p}, \qd \bz= \ll(\sum_{i=y}^N \hat{v}_i \rr)^{1- p}, \qd r=\frac{q}{p} \geq 1, 
$$
and then $2^{1/q- 1/p} B^*\leq B_*$.
So far, the proof of the DD-case in Theorem \ref{main thm} and the part (3) of Remark \ref{main rmk} is completed. 
\deprf

\section{Proof of the NN-case.}

In this section, we consider the NN boundary condition of weighted Hardy inequalities. Before further analysis, we give some properties about the constant $m(\mathbf{x})$ of a given sequence $\mathbf{x}$, which is defined by (\ref{NN boundary}). The first result is the existence and uniqueness.

\prp\lb{unique}
{\cms
For any sequence $\mathbf{x}$ on $[-M, N]$, there exists a unique constant $m= m(\mathbf{x})$ such that
$$
\sum_{n= -M}^N u_n \ll| x_n- m \rr|^{q-2} \ll(x_n- m \rr) =0.
$$
Moreover, the constant $m$ satisfies
\be\lb{inf}
\sum_{n= -M}^N u_n \ll|x_n -m \rr|^q = \inf_{c\in \mathbb{R}} \sum_{n= -M}^N u_n \ll|x_n -c \rr|^q. 
\de
}
\deprp

\prf
(a) For the existence, define a continuous function 
\be\lb{f}
f(t) = \sum_{n= -M}^N u_n \ll| x_n- t \rr|^{q-2} \ll(x_n- t \rr), 
\de
then we have $f(\max(\mathbf{x} ))\leq 0$ and $f(\min(\mathbf{x} ))\geq 0$, where $\max(\mathbf{x}):= \max \ll\{ x_n: -M\leq n\leq N \rr\}$ and $\min(\mathbf{x}):= \min \ll\{x_n: -M\leq n\leq N \rr\}$. Therefore, there exists $\min(\mathbf{x})\leq m\leq \max(\mathbf{x})$ such that $f(m)= 0$. 

(b) For the uniqueness, we prove it by contradiction. Assume that there are two constants $m_1< m_2$ satisfying $f(m_1)= f(m_2)= 0$. Write
$$
\aligned
A_1 &= \ll\{ n\in [-M, N]: \min(\mathbf{x})\leq x_n \leq m_1 \rr\}, \\
A_2 &= \ll\{ n\in [-M, N]: m_1 < x_n < m_2 \rr\}, \\
A_3 &= \ll\{ n\in [-M, N]: m_2 \leq x_n \leq \max(\mathbf{x}) \rr\},
\endaligned
$$
then we have 
$$
\aligned
0 &= f(m_1) \\
&= (- 1)\sum_{n\in A_1} u_n \ll|x_n -m_1 \rr|^{q-1} + \sum_{n\in A_2} u_n \ll|x_n -m_1 \rr|^{q-1}  \\
& \qqd + \sum_{n\in A_3} u_n \ll|x_n -m_1 \rr|^{q-1} \\
&> (- 1) \sum_{n\in A_1} u_n \ll|x_n -m_2 \rr|^{q-1}+ (- 1) \sum_{n\in A_2} u_n \ll|x_n -m_2 \rr|^{q-1} \\
& \qqd + \sum_{n\in A_3} u_n \ll|x_n -m_2 \rr|^{q-1} \\
&= f(m_2)=0.
\endaligned
$$
This leads to a contradiction. 

(c) We are now turning to prove (\ref{inf}). Define a continuous function
\be\lb{F}
F(t)= \sum_{n= -M}^N u_n |x_n -t|^q. 
\de
Taking the derivative with respect to $t$, the derivative of absolute value function means $(|x|)' = \sgn (x)$. Then we have $F'(t)= (-q) f(t)$ where $f$ is defined by (\ref{f}). By part (a) and (b), there is an unique constant $m$ such that $f(m)=0$. Using similar methods of part (b), for any $\vz>0$, we have
$$
f(m- \vz)> 0 \qd \text{and} \qd f(m+ \vz)< 0.
$$
It means that $m$ is an extreme point, moreover, $F$ reaches the minimum at $m$.
\deprf

To study (\ref{NN Hardy}), we start with the splitting technique again. To save our notations, without any confusion, we use the similar notations of Section 3 with different meanings. For any $\uz \in [-M, N]$ and $\gz \in [0, 1]$, define $u_n^-$, $v_n^-$ on $[-M+1, \uz]$ and $u_n^+$, $v_n^+$ on $[\uz, N]$ as 
\be\lb{uv-}
u_n^-= u_{n- 1}, n\in [-M+1, \uz], \qqd u_n^+= u_n, n\in [\uz, N],
\de
and
\be\lb{uv+}
v_n^- = 
\begin{cases}
v_n, & n\in [-M+1, \uz-1], \\
\gz^{1-p} v_\uz, & n=\uz,
\end{cases}
v_n^+ = 
\begin{cases}
(1- \gz)^{1-p} v_\uz, & n=\uz, \\
v_n, & n\in [\uz+1, N]. 
\end{cases}
\de
Before using the splitting technique, we give the following result.  
\lmm\lb{NN x+-}
{\cms
For any $\uz\in [-M, N]$, $\gz\in [0, 1]$ and sequence $\mathbf{x}$, define $\mathbf{x}^-$ on $[-M+ 1, \uz+ 1]$ and $\mathbf{x}^+$ on $[\uz- 1, N]$ as 
$$
\aligned
x_n^- &= 
\begin{cases}
x_{n-1}, & n\in [-M+1, \uz], \\
(1- \gz) x_{\uz- 1} + \gz x_{\uz}, & n= \uz+1, 
\end{cases} \\
x_n^+ &= 
\begin{cases}
x_{n}, & n\in [\uz, N], \\
(1- \gz) x_{\uz- 1} + \gz x_{\uz}, & n= \uz-1.  
\end{cases}
\endaligned
$$
Then we have 
$$
\aligned
\sum_{n=-M}^N u_n |x_n - m(\mathbf{x})|^q &= \sum_{n=-M+ 1}^\uz u_n^- |x_n^- - m(\mathbf{x})|^q + \sum_{n=\uz}^N u_n^+ |x_n^+ - m(\mathbf{x})|^q, \\
\sum_{n= -M}^N v_n |x_n - x_{n-1} |^p &= \sum_{n= -M+1}^\uz v_n^- |x_n^- - x_{n+1}^- |^p + \sum_{n= \uz}^N v_n^+ |x_n^+ - x_{n-1}^+ |^p. 
\endaligned
$$
}
\delmm

Here are the inequalities with single boundary condition. 
\begin{align}
\ll[\sum_{n= -M+ 1}^{\uz} u_n^- |x_n|^q \rr]^{1/q}\leq &A^-(\uz, \gz) \ll[\sum_{n= -M+1}^{\uz} v_n^- |x_n - x_{n+ 1} |^p \rr]^{1/p}, \  x_{\uz+ 1}= 0, \lb{NN A-}\\
\ll[\sum_{n= \uz}^{N} u_n^+ |x_n|^q \rr]^{1/q}\leq &A^+(\uz, \gz) \ll[\sum_{n= \uz}^N v_n^+ |x_n - x_{n- 1} |^p \rr]^{1/p}, \  x_{\uz}= 0, \lb{NN A+}. 
\end{align}
Without any confusion, we use the same notations $A$ and $A_{\uz}^\pm$ to express the optimal constants of the corresponding inequalities. Similarly to Proposition \ref{relationship}, the optimal constant $A$ is controlled by $A_{\uz}^\pm$. 

\prp\lb{NN relationship}
{\cms
For $1\leq p\leq q < \infty$, we have 
$$
\aligned
2^{1/q- 1/p} &\sup_{\uz\in [-M, N]} \sup_{\gz\in [0, 1]} \ll(A^- (\uz, \gz) \wedge A^+ (\uz, \gz) \rr) \leq A  \\ 
&\leq \inf_{\uz\in [-M, N]} \inf_{\gz\in [0, 1]} \ll(A^- (\uz, \gz) \vee A^+ (\uz, \gz) \rr).
\endaligned
$$
}
\deprp

\prf (a) For any $\uz \in [-M, N]$, $\gz\in [0, 1]$ and sequence $\mathbf{x}$, set 
$$
\tilde{x}_n = x_n - \ll( (1- \gz) x_\uz + \gz x_{\uz+ 1} \rr), -M \leq n \leq N.  
$$ 
By Lemma \ref{NN x+-}, $c_r$-inequality and Proposition \ref{unique}, we obtain
$$
\aligned
\sum_{n= -M}^N& v_n |x_n- x_{n-1}|^p = \sum_{n= -M+1}^\uz v_n^- |\tilde{x}_n^- - \tilde{x}_{n+1}^- |^p + \sum_{n= \uz}^N v_n^+ |\tilde{x}_n^+ - \tilde{x}_{n-1}^+ |^p \\
&\geq \ll(2^{(p/q -1) \vee 0} \rr)^{-1} \ll[\ll(A^-(\uz, \gz) \rr)^{-p} \wedge \ll(A^+(\uz, \gz) \rr)^{-p} \rr] \ll(\sum_{n= -M}^N u_n |\tilde{x}_n|^q \rr)^{p/q} \\
&\geq \ll[\ll(A^-(\uz, \gz) \rr)^{-p} \wedge \ll(A^+(\uz, \gz) \rr)^{-p} \rr] \ll(\sum_{n= -M}^N u_n |x_n - m(\mathbf{x}) |^q \rr)^{p/q}. 
\endaligned
$$
Since $\uz$, $\gz$ and $\mathbf{x}$ are arbitrary, we obtain the upper bound of $A$
$$
A \leq \inf_{\uz\in [-M, N]} \inf_{\gz\in [0, 1]} \ll(A^- (\uz, \gz) \vee A^+ (\uz, \gz) \rr).
$$

(b) Fix $\uz \in [-M, N]$ and $\gz \in [0, 1]$ and $\vz>0$. Let $\mathbf{x}^{(1)}$ be a positive sequence on $[-M+ 1, \uz+ 1]$ satisfies 
$$
\aligned
&\sum_{n= -M+1}^\uz u_n^- \ll| x_n^{(1)} \rr|^q = 1, \qd x_{\uz+ 1}^{(1)}= 0, \\
& \ll[ \sum_{n=-M+ 1}^\uz v_n^-  \ll| x_n^{(1)}- x_{n+1}^{(1)} \rr|^p \rr]^{1/p} < \ll(A^-(\uz, \gz)\rr)^{-1} + \vz. 
\endaligned
$$
Let $\mathbf{x}^{(2)}$ be a positive sequence on $[\uz- 1, N]$ satisfies 
$$
\aligned
& \sum_{n= \uz}^N u_n^+ \ll| x_n^{(2)} \rr|^q = 1, \qd x_{\uz- 1}^{(2)}= 0, \\
& \ll[  \sum_{n=\uz}^N v_n^+ |x_n^{(2)}- x_{n- 1}^{(2)}|^p \rr]^{1/p} < \ll(A^+(\uz, \gz)\rr)^{-1} + \vz. 
\endaligned
$$
Here we define sequence $\mathbf{x}$ on $[-M, N]$ by $\mathbf{x}^{(1)}$ and $\mathbf{x}^{(2)}$. Set
$$
x_n= 
\begin{cases}
-c x_{n+1}^{(1)}, & -M\leq n\leq \uz- 1, \\
x_n^{(2)}, & \uz\leq n\leq N, 
\end{cases}
$$
where 
$$
c= \ll( \frac{\sum_{n= \uz}^N u_n^+ |x_n^{(2)}|^{q- 1} }{\sum_{n= -M}^\uz u_n^- |x_n^{(1)}|^{q- 1}} \rr)^{1/(q- 1)}. 
$$
Obviously, $\mathbf{x}$ satisfies $\sum_{n= -M}^N u_n |x_n|^{q- 2} (x_n) = 0$, by Proposition \ref{unique} we have $m(\mathbf{x})= 0$. Moreover, we have
$$
\sum_{n= -M}^N u_n |x_n|^{q} = \sum_{n= -M+ 1}^\uz u_n^- c^q |x_n^{(1)}|^{q} + \sum_{n= \uz}^N u_n^+ |x_n^{(2)}|^{q} = c^q +1. 
$$

Note that the function 
$$
F(\gz)= \gz^{1- p} a^p+ (1- \gz)^{1- p} b^p, \qd \gz \in [0, 1], a>0, b>0, p\in (1, \infty)
$$
achieves its minimum $(a+ b)^p$ at $\gz^*= \dfrac{a}{a+ b}$. Applying this result with $a= c x_\uz^{(1)}$ and $b= x_\uz^{(2)}$, we obtain 
\be\lb{NN =1}
v_\uz |x_\uz^{(2)}+ c x_\uz^{(1)}|^p \leq \gz^{1- p} v_\uz |c x_\uz^{(1)}|^p + (1- \gz)^{1- p} v_\uz |x_\uz^{(2)}|^p. 
\de
Hence, 
$$
\aligned
\sum_{n= -M}^N v_n |x_n - x_{n- 1}|^p &\leq \sum_{n= -M}^{\uz} v_n^- c^p |x_n^{(1)} - x_{n+ 1}^{(1)} |^p + \sum_{n= \uz}^N v_n^+ |x_n^{(2)} - x_{n- 1}^{(2)} |^p \\
&\leq 2^{\frac{1}{p}- \frac{1}{q}} \ll(c^q +1 \rr)^{\frac{1}{q}} \ll[(A^- (\uz, \gz) )^{-1} \vee (A^+ (\uz, \gz) )^{-1}+ \vz \rr].
\endaligned
$$
Since $\uz$, $\gz$ and $\vz$ are arbitrary, we obtain the lower bound of $A$
$$
A \geq 2^{1/q- 1/p} \sup_{\uz\in [-M, N]} \sup_{\gz\in [0, 1]} \ll( A^- (\uz, \gz) \wedge A^+ (\uz, \gz) \rr). 
$$
\deprf

Using the results of ND-case and the DN-case of Theorem \ref{main thm}. Hence, we have the basic estimate of $A^{\pm}$  
\be\lb{NN estimate}
B^\pm (\uz, \gz)\leq A^\pm (\uz, \gz)\leq k_{q, p} B^\pm (\uz, \gz), 
\de
where $k_{q, p}$ is defined as (\ref{kqp}) and $B^\pm (\uz, \gz)$ are
\begin{align}\lb{NN B^+-}
B^- (\uz, \gz) &= \sup_{n\in [-M, \uz- 1]} \ll(\sum_{i= -M}^n u_i \rr)^{1/q} \ll(\sum_{j= n+1}^{\uz- 1} \hat{v}_j + \gz \hat{v}_\uz \rr)^{1/p^*}, \\
B^+ (\uz, \gz) &= \sup_{n\in [\uz, N]} \ll(\sum_{j= \uz+ 1}^n \hat{v}_j + (1- \gz) \hat{v}_\uz \rr)^{1/p^*} \ll(\sum_{i= n}^{N} u_i \rr)^{1/q}.
\end{align}

Having these preparations at hand, we are ready to give the basic estimate of the NN-case. 

\medskip
\noindent {\bf Proof of the NN-case.} For brevity, we use $B^*$ (resp., $B_*$) instead of $B^{* NN}$ (resp., $B_*^NN$) in this proof. 

(a) An argument similar to part (a) and (b) in the proof of the DD-case, we have 
$$
\sup_{\uz \in [-M, N]} \sup_{0\leq \gz \leq 1} \ll[B^- (\uz, \gz) \wedge B^+ (\uz, \gz) \rr] \leq B^*, 
$$
and there exist $(\bar{\uz}, \bar{\gz})$ such that 
$$
B^- (\bar{\uz}, \bar{\gz})= B^+ (\bar{\uz}, \bar{\gz}). 
$$
Furthermore, by Proposition \ref{NN relationship} and (\ref{NN estimate}), we obtain 
$$
\aligned
A &\leq \inf_{\uz\in [-M, N]} \inf_{0\leq \gz\leq 1} \ll(A^-(\uz, \gz) \vee A^+(\uz,\gz) \rr) \\
&\leq k_{q,p} \inf_{\uz\in [-M, N]} \inf_{0\leq \gz\leq 1} \ll(B^-(\uz, \gz) \vee B^+(\uz,\gz) \rr) \\
&\leq k_{q,p} B^-(\bar{\uz}, \bar{\gz}) \leq k_{q,p} B^*.
\endaligned
$$
This method is quite the same as the proof (c) of DD-case. That is the upper bound $A\leq k_{q,p} B^*$. 

(b) The next goal is to estimate the lower bound of $A$. Fix $x, y\in [-M, N]$ with $x< y$. We assert that exist $\bar{\uz} \in [x, y]$ and $\bar{\gz} \in [0, 1]$ such that 
\be\lb{NN =3}
\sum_{n= -M}^{\bar{\uz} - 1} u_n \ll|\sum_{i= n \vee x}^{\bar{\uz}- 1} \widehat{v_{i+1}^-} \rr|^{q- 1} = \sum_{n= \bar{\uz}}^N u_n \ll|\sum_{i= \bar{\uz}}^{y \wedge n} \widehat{v_i^+} \rr|^{q- 1}.
\de
In fact, the idea is similar to the proof (b) of DD-case. Let 
$$
C^- (\uz, \gz) = \sum_{n= -M}^{\uz - 1} u_n \ll|\sum_{i= n \vee x}^{\uz- 1} \widehat{v_{i+1}^-} \rr|^{q- 1}, \qd C^+ (\uz, \gz) = \sum_{n= \uz}^N u_n \ll|\sum_{i= \uz}^{y \wedge n} \widehat{v_i^+} \rr|^{q- 1}. 
$$
By direct calculation, we have $C^-(\uz, 0)= C^-(\uz-1, 1)$ and $C^+(\uz, 0)= C^+(\uz- 1, 1)$. Besides, for any $\gz\in [0, 1]$, when $\uz$ varies from $x$ to $y$, $C^- (\uz, \gz)$ goes from $0$ to a positive number and $C^+ (\uz, \gz)$ goes from a positive number to $0$. Hence, let 
$$
\bar{\uz}+ 1= \min \{x\leq n\leq y: C^-(n, 0)\geq C^+(n, 0) \}, 
$$
then we have
\be\lb{NN =4}
\begin{cases}
C^-(\bar{\uz}, 0)\leq C^+(\bar{\uz}, 0), \\
C^-(\bar{\uz}, 1)\geq C^+(\bar{\uz}, 1).
\end{cases}
\de
Hence, the existence of $(\bar{\uz}, \bar{\gz})$ is clear because when $\gz$ varies from $0$ to $1$, $C^- (\bar{\uz}, \gz)$ increases continuously, correspondingly, $C^+ (\bar{\uz}, \gz)$ decreases continuously. 

Fix this $(\bar{\uz}, \bar{\gz})$, define $\mathbf{\bar{x}}$ on $[-M, N]$ 
\be\lb{bar x}
\bar{x}_n =
\begin{cases}
(-1) \sum_{i= n \vee x}^{\bar{\uz} - 1} \widehat{v_{i+ 1}^-}, & -M\leq n\leq \bar{\uz}- 1, \\
\sum_{i= \bar{\uz}}^{y \wedge n} \widehat{v_i^+}, & \bar{\uz} \leq n \leq N.
\end{cases}
\de
Obviously, we have $m(\mathbf{\bar{x}})= 0$. 

By direct calculation, we obtain 
$$
\sum_{n=-M}^N v_n \ll| \bar{x}_n - \bar{x}_{n- 1} \rr|^p = \sum_{n= x+ 1}^y \hat{v}_n, 
$$
and 
$$
\sum_{n= -M}^N u_n \ll|\bar{x}_n - m(\mathbf{\bar{x}}) \rr|^q \geq \ll(\sum_{i= -M}^x u_i \rr) \ll(\sum_{j=x+ 1}^{\bar{\uz}} \widehat{v_j^-} \rr)^q + \ll(\sum_{i= y}^N u_i \rr) \ll(\sum_{j=\bar{\uz}}^y \widehat{v_j^+} \rr)^q. 
$$
Hence, we have
\be\lb{NN =2} 
A \geq \ll(\sum_{n= x+1}^y \hat{v}_n \rr)^{-1/p} \ll[ \ll(\sum_{i= -M}^x u_i \rr) \ll(\sum_{i=x+ 1}^{\bar{\uz}} \widehat{v_i^-} \rr)^q + \ll(\sum_{i= y}^N u_i \rr) \ll(\sum_{i=\bar{\uz}}^y \widehat{v_i^+} \rr)^q \rr]^{1/q}. 
\de
Note that the function $F(x)= \az x^q + \bz (1- x)^q$ achieves its minimum $ \ll(\az^{1- q^*} + \bz^{1- q^*} \rr)^{1- q}$ at $x_0= \bz^{q^*- 1} \ll(\az^{q^*- 1}+ \bz^{q^* -1} \rr)^{-1}$. Using this result with 
$$
\az= \sum_{i= -M}^x u_i, \qd \bz= \sum_{i= y}^N u_i,
$$
we have 
$$
A \geq \ll(\sum_{n= x+1}^y \hat{v}_i \rr)^{1/p^*} \ll[ \ll(\sum_{i=-M}^x u_i \rr)^{1- q^*} + \ll(\sum_{i= y}^N u_i \rr)^{1- q^*} \rr]^{-1/q^*}. 
$$
Since $x$ and $y$ are arbitrary, we obtain the lower bound of $A$.

(c) For the relationship of $B^*$ and $B_*$, we can work it out in the same way of part (d) in the proof (d) of DD-case. 
\deprf

\section{Example}

As an application of Theorem \ref{main thm}, the first example is the dual version of \rf{Liao}{Example 5.2}. The optimal constant is explicit in this example.  
\xmp\lb{example1}
Let $1< p\leq q< \infty$ and $N= \infty$. For $n\geq 1$, define
$$
u_n \equiv 1, \qd v_n = \ll[ n^{-p^*/ q} - (n+1)^{-p^*/q} \rr]^{1- p}. 
$$
Then the optimal constant of {\rm ({\ref{ND Hardy}})} is $A^{ND}= k_{q, p}$, and the basic estimate is
$$
B^{ND}=1, \qd 1 \leq A^{ND} = k_{q, p}. 
$$
\dexmp

The next example is modified from \rf{Opic}{Example 8.6}. We will use the DD-case of Theorem \ref{main thm} to give the necessary and sufficient conditions for the validity of the weighted Hardy inequality. 

\xmp\lb{example2}
Let $1< p\leq q < \infty$, $\az, \bz \in \mathbb{R}^+$. We consider the following Hardy inequality
\be\lb{ex2}
\ll[\sum_{n= 1}^\infty \frac{1}{n^\az} |x_n|^q \rr]^{1/q} \leq A \ll[\sum_{n=1}^\infty n^\bz |x_n- x_{n- 1}|^p \rr]^{1/p},
\de
where $x_{\infty} = x_0 = 0$.  \\
{\rm (1)} If $\bz= p- 1$, then this inequality holds if and only if $\az > 1$.  \\
{\rm (2)} If $\bz \neq p- 1$, then this inequality holds if and only if $\az \geq 1+ \dfrac{q}{p} (p- 1- \bz) $. 
\dexmp

\prf
(1) When $0< \bz < p-1$, we have $-1 < \bz (1- p^*)< 0$, which means the sequence $\sum_{n= y}^\infty n^{\bz (1- p^*)}$ is divergent. By the definition of $B^{* DD}$ and $B_*^{DD}$, we obtain
$$
\aligned
B^{* DD}= B_*^{DD} &= \sup_{1\leq x< y< \infty} \ll(\sum_{n= x}^{y- 1} \frac{1}{n^\az} \rr)^{1/q} \ll(\sum_{n= 1}^x n^{\bz (1- p^*) } \rr)^{1/p^*} \\
&\sim x^{\frac{1- \az}{q}+ \frac{\bz(1- p^*)+ 1}{p^*}}, \qd \text{as $x\rightarrow \infty$}. 
\endaligned
$$
Hence, the necessary and sufficient conditions of $B^{* DD}< \infty$ are
$$
\az > 1 \qd \text{and} \qd \frac{1- \az}{q}+ \frac{\bz (1- p^*)+ 1}{p^*} \leq 0. 
$$
That is $\az \geq 1+ \dfrac{q}{p} (p- 1- \bz) $, since $p-1 -\bz > 0$. 

(2) When $\bz= p- 1$, we have $\bz (1- p^*)= -1$. An argument similar to part (1), we have
$$
B^{* DD}= B_*^{DD} \sim x^{\frac{1- \az}{q}} (\ln x)^{1/{p*}}, \qd \text{as $x\rightarrow \infty$}. 
$$
Obviously, the necessary and sufficient condition is $\az > 1$. 

(3) Now lets consider the case $\bz > p-1$. In this case, we have $\bz (1- p^*)< -1$. By the definition, we have  
$$
\aligned
(B^{* DD})^{-q} \geq \Bigg\{ & \sup_{1\leq x < y< \infty} \ll(\sum_{n= x}^{y- 1} \frac{1}{n^\az} \rr)^{-1} \ll(\sum_{n= 1}^x n^{\bz (1- p^*)} \rr)^{-q/p^*} \Bigg\}^{-1} \\
& + \Bigg\{ \sup_{1\leq x < y< \infty} \ll(\sum_{n= x}^{y- 1} \frac{1}{n^\az} \rr)^{-1} \ll(\sum_{n= y}^\infty n^{\bz (1- p^*)} \rr)^{-q/p^*} \Bigg\}^{-1}, 
\endaligned
$$
hence, the sufficient condition of $B^{* DD} < \infty$ is
\begin{align}
\sup_{1\leq x < y< \infty} & \ll( \sum_{n= x}^{y- 1} \frac{1}{n^\az} \rr) \ll(\sum_{n= 1}^x n^{\bz (1- p^*)} \rr)^{q/p^*} < \infty, \text{\qd or \qd} \lb{ex2=1} \\
\sup_{1\leq x < y< \infty} & \ll(\sum_{n= x}^{y- 1} \frac{1}{n^\az} \rr) \ll(\sum_{n= y}^\infty n^{\bz (1- p^*)} \rr)^{q/p^*} < \infty. \lb{ex2=2} 
\end{align}
For the first part (\ref{ex2=1}), we have
$$
\sup_{1\leq x < y< \infty}  \ll( \sum_{n= x}^{y- 1} \frac{1}{n^\az} \rr) \ll(\sum_{n= 1}^x n^{\bz (1- p^*)} \rr)^{q/p^*} \sim x^{1- \az+ \frac{q}{p^*} - \frac{\bz q}{p}} \qd \text{as $x\rightarrow \infty$}.  
$$
Hence, we need 
$$
\az >1 \text{\qd and \qd} \az \geq 1+ \dfrac{q}{p} (p- 1- \bz).
$$
Since $p- 1- \bz < 0$, we obtain $\az > 1$. For the second part (\ref{ex2=2}), by the similar approach, we have
$$
\sup_{1\leq x < y< \infty} \ll(\sum_{n= x}^{y- 1} \frac{1}{n^\az} \rr)  \ll(\sum_{n= y}^\infty n^{\bz (1- p^*)} \rr)^{q/p^*} \sim y^{1- \az+ \frac{q}{p^*} - \frac{\bz q}{p}} \qd \text{as $y\rightarrow \infty$}. 
$$
Hence, we need
$$
\az \geq 1+ \dfrac{q}{p} (p- 1- \bz).
$$
In conclusion of the first and the second part, $\az \geq 1+ \frac{q}{p^*}- \frac{\bz q}{p}$ implies $B^{* DD} < \infty$ when $\bz > p-1$. By Remark \ref{main rmk}, we have $B_*^{DD} \leq B^{* DD} < \infty$. 

And now we consider the necessary condition of $B^{* DD} < \infty$. We have 
$$
\aligned
\infty > B^{* DD} & \geq \sup_{1< y< \infty} \ll(\sum_{n= 1}^{y- 1} \frac{1}{n^\az} \rr)^{1/q} \ll[1 + \ll(\sum_{n=y}^\infty n^{\bz (1- p^*)} \rr)^{-q/{p^*}} \rr]^{-1/q} \\
&\sim y^{\frac{1- \az}{q} + \frac{1}{p^*}(\bz (1- p^*) +1)} \qd \text{as $y\rightarrow \infty$}. 
\endaligned
$$
The necessary condition such that the above inequality holds is  
$$
\az \geq 1+ \frac{q}{p^*}- \frac{\bz q}{p}. 
$$

Hence, when $\bz> p- 1$, the necessary and sufficient condition for the validity of this weighted Hardy inequality is $\az \geq 1+ \frac{q}{p^*}- \frac{\bz q}{p}$.
\deprf

\rmk
{\cms
To compare the results for the discrete case with that of continuous one, we excerpt the example of \rf{Opic}{Example 8.6}: the Hardy inequality
$$
\ll(\int_{0}^{\infty} \frac{1}{x^\az} |f(x)|^q \d x \rr)^{1/q} \leq A \ll(\int_{0}^{\infty} x^{\bz} |f'(x)|^p \d x \rr)^{1/p}, f(0)= \lim_{t \rightarrow \infty} f(t)= 0
$$
holds for every absolutely continuous function with compact support with a finite constant $A$ if and only if 
$$
\bz \neq p- 1, \qqd \az = 1- \frac{\bz q}{p} + \frac{q}{p^*}.
$$
}
\dermk

As we know, weighted Hardy inequalities' applications have been expanded to probability theory, cf. \cite{Chen5}. The following example comes from the birth-death processes with reflecting boundaries at origin and infinity, cf. \rf{Chen5}{Example 6.7}. The basic estimates of this inequality will be presented by Theorem \ref{main thm}. 

\xmp\lb{example3}
Let $1< p\leq q < \infty$, $N< \infty$. We consider the following inequality
\be\lb{ex3}
\ll[\sum_{n= 1}^N r^n |x_n- m(\mathbf{x}) |^q \rr]^{1/q} \leq A \ll[\sum_{n=1}^\infty b r^n |x_n- x_{n- 1}|^p \rr]^{1/p},
\de
where $x_0= x_1$, $b> 0$, $0< r< 1$ and $m(\mathbf{x})$ is the constant defined by (\ref{NN boundary}). Then we have
$$
B_*^{NN} \leq A \leq k_{q, p} B^{* NN}, 
$$
where
$$
\aligned
B^{* NN} &= b^{-1/p} \ll(r^{1- p^*} - 1 \rr)^{-1/p^*} \ll[\frac{r^{N (1- p^*)}- r^{2 (1- p^*)}}{r^{-p^*/q}+ r^{N (-p^*/q)}} \rr]^{1/p^*}, \\
B_*^{NN} &= b^{-1/p} \ll(r^{1- p^*} - 1 \rr)^{-1/p^*} \ll[\frac{\ll( r^{N (1- p^*)}- r^{2 (1- p^*)}\rr)^{q^*/p^*}}{r^{1- q^*}+ r^{N (1- q^*)}} \rr]^{1/q^*}. 
\endaligned
$$
Obviously, when $p= q= 2$ and $N= \infty$, this result is consistent with the Example 6.7 in \cite{Chen5}. 
\dexmp

\prf
We could compute it directly by Theorem \ref{main thm}, we have
$$
\aligned
B^{* NN} & = \sup_{1\leq x< y \leq N} \ll(\sum_{n= x+ 1}^y b^{1- p^*} r^{n(1- p^*)} \rr)^{1/p^*} \ll[\ll(\sum_{n= 1}^x r^n \rr)^{-p^*/q} + \ll(\sum_{n= y}^N r^n \rr)^{-p^*/q} \rr]^{-1/p^*} \\
&= \sup_{1\leq x< y \leq N} b^{-1/p} \ll(r^{1- p^*} - 1 \rr)^{-1/p^*} \ll(1- r\rr)^{-1/q} F(x, y)^{1/p^*}. 
\endaligned
$$
where
$$
F(x, y) := \frac{r^{y (1- p^*)}- r^{(x+ 1) (1- p^*)}}{(r- r^{x+ 1})^{-p^*/q}+ (r^y- r^{N+ 1})^{-p^*/q}}.
$$
Since $1< p\leq q< \infty$, by the derivation of $F(x, y)$, it is easy to show the $F(x, y)$ is decreasing with respect to $x$ and is increasing with respect to $y$. That means $\sup\limits_{1\leq x< y \leq N} F(x, y)= F(1, N)$. 

For the lower estimates of $A$, we have 
$$
\aligned
B_*^{NN} & = \sup_{1\leq x< y \leq N} \ll(\sum_{n= x+ 1}^y b^{1- p^*} r^{n(1- p^*)} \rr)^{1/p^*} \ll[\ll(\sum_{n= 1}^x r^n \rr)^{1- q^*} + \ll(\sum_{n= y}^N r^n \rr)^{1- q^*} \rr]^{-1/p^*} \\
&= \sup_{1\leq x< y \leq N} b^{-1/p} \ll(r^{1- p^*} - 1 \rr)^{-1/p^*} \ll(1- r\rr)^{-1/q} F_0 (x, y)^{1/q^*}. 
\endaligned
$$
where
$$
F_0 (x, y) := \frac{\ll[ r^{y (1- p^*)}- r^{(x+ 1) (1- p^*)} \rr]^{q^*/p^*}}{(r- r^{x+ 1})^{1- q^*}+ (r^y- r^{N+ 1})^{1- q^*}}.
$$
In a similar way, $F_0 (x, y)$ is decreasing with respect to $x$ and is increasing with respect to $y$. Hence, we obtain $\sup\limits_{1\leq x< y \leq N} F_0 (x, y)= F_0 (1, N)$. 

\deprf

\noindent {\bf Acknowledgements} $\qd$ This paper is based on the series of studies of my supervisor Prof. M. F. Chen. Heartfelt thanks are given to my supervisor for his careful guidance and helpful suggestions. Thanks are also given to Prof. Y. H. Mao, Prof. F. Y. Wang and Prof. Y. H. Zhang for their comments and suggestions, which lead to lots of improvements of this paper. 

The research is supported by NSFC (Grant No. 11131003) and by the ``985'' project from the Ministry of Education in China.


\begin{thebibliography}{10}
\setlength{\itemsep}{-0.8ex}

\bibitem{Bennett3} G. Bennett,
{\it Some elementary inequalities III}, 
Quart. J. Math. Oxford (2). {\bf 42} (1991), 149--174.
\lb{Bennett3}

{\small
\bibitem{Bliss}
G. A. Bliss, 
An integral inequality.
Journ. L.M.S. 5, 40-46 (1930)
\lb{Bliss}

\bibitem{Buttazzo} G. Buttazzo, M. Giaquinta, S. Hildebrandt, 
One-dimensional variational problems. An introduction, 
Oxford Lecture Series in Mathematics and its Applications, 15, 
The Clarendon Press, Oxford University Press, New York, 1998.
\lb{Buttazzo}

\bibitem{Chen1}
M. F. Chen,
Explicit bounds of the first eigenvalue. 
Sci. China {\rm (A)} 43(10), 1051--1059 (2000)
\lb{Chen1}

\bibitem{Chen6}
M. F. Chen,
Eigenvalues, Inequalities and Ergodic Theory, 
Springer, London, 2005
\lb{Chen6}

\bibitem{Chen5}
M. F. Chen,
Speed of stability for birth-death processes. 
Front. Math. China, 5(3), 379-515 (2010)
\lb{Chen5}

\bibitem{Chen9}
M. F. Chen,
Basic estimates of stability rate for one-dimensional diffusions. 
Chapter 6 in "Probability Approximations and Beyond", eds. A. D. Barbour, H. P. Chan and D. Siegmund, Lecture Notes in Statistics, 205, (2012) 75-99
\lb{Chen9}

\bibitem{Chen4}
M. F. Chen,
Bilateral Hardy-type inequalities.
Acta Math. Sin. Eng. Ser. 29:1, 1-32 (2013)
\lb{Chen4}

\bibitem{Chen3}
M. F. Chen,
The optimal constant in Hardy-type Inequalities.
preprint, (2014)
\lb{Chen3}

\bibitem{Chen7}
M. F. Chen, F. Y. Wang,
Cheeger's inequalities for general symmetric forms and existence criteria for spectral gap.  
Abstract. Chin. Sci. Bulletin 43:18, 1516-1519. Ann. Prob., 28:1, 235-257 (2000)
\lb{Chen7}

\bibitem{Chen2}
M. F. Chen, L. D. Wang, Y. H. Zhang,
Mixed principal eigenvalues in dimension one.
Frontier of Mathematic in China, 8(2), 317--343 (2013)
\lb{Chen2}

\bibitem{Chen8}
M. F. Chen, Y. H. Zhang, X. L. Zhao,
Dual variational formulas for the first Dirichlet eigenvalue on half-line.
Sci. China 46:6, 847-861. (2003)
\lb{Chen8}

\bibitem{Gurka} 
P. Gurka,  
Generalized Hardy's inequality,  
\v{C}asopis P\v{e}st. Mat. {\bf 109} (1984), no. 2, 194--203. 
\lb{Gurka}

\bibitem{Wang}
L. D. Wang, 
Mixed eigenvalues for $p$-Laplacian and converge speed for several kinds of Markov processes.
Ph.D. Thesis, Beijing Normal University, Beijing, China, 2013
\lb{Wang}

\bibitem{Hardy_book}
G. H. Hardy, J. E. Littlewood, G. Polya,
Inequalities, 
2nd edition, Cambridge University Press, 1967.
\lb{Hardy_book}

\bibitem{Kufner}
A. Kufner, L. Maligranda, L. E. Persson, 
The Hardy Inequality: About its History and Some Related Results, 
Vydavatelsky Servis, 2007 
\lb{Kufner}

\bibitem{Liao}
Z. W. Liao,
Discrete Hardy-type Inequalities. 
(2014)
arXiv: 1406.1984v2 [math. FA] 
\lb{Liao}

\bibitem{Manakov}
V. M. Manakov,
On the best constant in weighted inequalities for Riemann-Liouville integrals.
Bull. London Math. Soc., 24, 442-448 (1992) 
\lb{Manakov}

\bibitem{Mao}
Y. H. Mao, 
Nash inequalities for Markov processes in dimension one. 
Acta Math. Sinica, 18(1): 147-156 (2002)
\lb{Mao}

\bibitem{Mao2}
Y. H. Mao, L. Y. Xia,
Spectral gap for jump processes by decomposition method.
Front. Math. China 4:2, 335-347 (2009)
\lb{Mao2}

\bibitem{Maz'ya}
V. Maz'ya,
Sobolev Space with Applications to Elliptic Partial Differential Equations (2nd Ed.),
Springer, Berlin, 2011
\lb{Maz'ya}

\bibitem{Opic}
B. Opic, A. Kufner,
Hardy-type Inequalities. 
Longman, New York, 1990
\lb{Opic}

}
\end{thebibliography}
\end{document}